\documentclass[11pt]{article}
\usepackage{amsfonts,latexsym,rawfonts,amsmath,amssymb,amsthm,graphicx}
\numberwithin{equation}{section}
\newtheorem{theorem}{Theorem}[section]
\newtheorem{lem}[theorem]{Lemma}
\newtheorem{thm}[theorem]{Theorem}
\newtheorem{pro}[theorem]{Proposition}
\newtheorem{cor}[theorem]{Corollary}

\newtheorem{rem}[theorem]{Remark}

\def\endproof{$\hfill\Box$\\}
\def\R{\mathbb{R}}

\title{A nonexistence theorem for proper biharmonic maps into general Riemannian manifolds}
\author{Volker Branding, Yong Luo}
\date{\today}

\begin{document}
\maketitle
\begin{abstract}
In this note we prove a nonexistence result for proper biharmonic maps from complete non-compact Riemannian manifolds
of dimension \(m=\dim M\geq 3\) with infinite volume that admit an Euclidean type Sobolev inequality
into general Riemannian manifolds by assuming finiteness of $\|\tau(\phi)\|_{L^p(M)}, p>1$ and smallness of $\|d\phi\|_{L^m(M)}$.
This is an improvement of a recent result of the first named author, where he assumed $2<p<m$.
As applications we also get several nonexistence results for proper biharmonic submersions from complete non-compact manifolds into  general Riemannian manifolds.

\end{abstract}
\section{Introduction}
Let $(M,g)$ be a Riemannian manifold and $(N,h)$ a Riemannian manifold without boundary. For a $W^{1,2}(M,N)$ map $\phi$, the energy density of $\phi$ is defined by
$$ e(\phi)=|d\phi|^2=\rm{Tr_g}(\phi^\ast h),$$
where $\phi^\ast h$ is the pullback of the metric tensor $h$. The energy functional of the map $\phi$ is defined as $$E(\phi)=\frac{1}{2}\int_Me(\phi)dv_g.$$
The Euler-Lagrange equation of $E(\phi)$ is $\tau(\phi)=\rm Tr_g\bar{\nabla} d\phi=0$ and $\tau(\phi)$ is called the \textbf{tension field} of $\phi$. A map is called a \textbf{harmonic map} if $\tau(\phi)=0$.

The theory of harmonic maps has many important applications in various fields of differential geometry, including minimal surface theory, complex geometry, see \cite{SY} for a survey.

Much effort has been paid in the last several decades to generalize the notion of harmonic maps. In 1983, Eells and Lemaire (\cite{EL}, see also \cite{ES}) proposed to consider the bienergy functional
$$E_2(\phi)=\frac{1}{2}\int_M|\tau(\phi)|^2dv_g$$ for smooth maps between Riemannian manifolds. Stationary points of the bienergy functional are called \textbf{biharmonic maps}.
We see that harmonic maps are biharmonic maps and even more, minimizers of the bienergy functional. In 1986, Jiang \cite{Ji} derived the first and second variational formulas of the bienergy functional and initiated the study of biharmonic maps. The Euler-Lagrange equation of $E_2(\phi)$ is given by
$$\tau_2(\phi):=-\Delta^\phi\tau(\phi)-\sum_{i=1}^mR^N(\tau(\phi), d\phi(e_i))d\phi(e_i)=0,$$
where $\Delta^\phi:=\sum_{i=1}^m(\bar{\nabla}_{e_i}\bar{\nabla}_{e_i}-\bar{\nabla}_{\nabla_{e_i}e_i})$. Here, $\nabla$ is the Levi-Civita connection on $(M,g)$, $\bar{\nabla}$ is the induced connection on the pullback bundle $\phi^{\ast}TN$, and $R^N$ is the Riemannian curvature tensor on $N$.

The first nonexistence result for biharmonic maps was obtained by Jiang \cite{Ji}. He proved that biharmonic maps from a compact, orientable Riemannian manifold into a Riemannian manifold of nonpositive curvature are harmonic.
Jiang's theorem is a direct application of the Weitzenb\"ock formula. If $\phi$ is biharmonic, then
\begin{eqnarray*}
-\frac{1}{2}\Delta|\tau(\phi)|^2&=&\langle-\Delta^\phi\tau(\phi), \tau(\phi)\rangle-|\bar{\nabla}\tau(\phi)|^2
\\&=&Tr_g\langle R^N(\tau(\phi), d\phi)d\phi, \tau(\phi)\rangle-|\bar{\nabla}\tau(\phi)|^2
\\&\leq&0.
\end{eqnarray*}
The maximum principle implies that $|\tau(\phi)|^2$ is constant. Therefore $\bar{\nabla}\tau(\phi)=0$ and so by
$$div\langle d\phi, \tau(\phi)\rangle=|\tau(\phi)|^2+\langle d\phi, \bar{\nabla}\tau(\phi)\rangle,$$
we deduce that $div\langle d\phi, \tau(\phi)\rangle=|\tau(\phi)|^2$. Then, by the divergence theorem, we have $\tau(\phi)=0$. Generalizations of this result by making use of similar ideas are given in \cite{On}.

If $M$ is non-compact, the maximum principle is no longer applicable. In this case, Baird et al. \cite{BFO} proved that biharmonic maps
from a complete non-compact Riemannian manifold with nonnegative Ricci curvature into a nonpositively curved manifold with finite bienergy are harmonic.
It is natural to ask whether we can abandon the curvature restriction on the domain manifold and weaken the integrability condition on the bienergy.
In this direction, Nakauchi et al. \cite{NUG} proved that biharmonic maps from a complete manifold to a nonpositively curved manifold are harmonic if ($p=2$)
\\(i) $\int_M|d\phi|^2dv_g<\infty$ and $\int_M|\tau(\phi)|^pdv_g<\infty$, or
\\(ii) $Vol(M, g)=\infty$ and $\int_M|\tau(\phi)|^pdv_g<\infty.$

Later Maeta \cite{Ma} generalized this result by assuming that $p\geq2$ and further generalizations are given by the second named author in \cite{Luo1}, \cite{Luo2}.

Recently, the first named author proved a nonexistence result for proper biharmonic maps from complete non-compact manifolds into general target manifolds \cite{Ba}, by only assuming that the sectional curvatures of the target manifold have an upper bound.
Explicitly, he proved the following theorem.
\begin{thm}[Branding]\label{Bra}
Suppose that $(M,g)$ is a complete non-compact Riemannian manifold of dimension \(m=\dim M\geq 3\)
whose Ricci curvature is bounded from below and with positive injectivity radius.
Let $\phi: (M,g)\to (N,h)$ be a smooth biharmonic map, where \(N\) is another Riemannian manifold.
Assume that the sectional curvatures of $N$ satisfy $K^N\leq A,$  where $A$ is a positive constant.
If $$\int_M|\tau(\phi)|^pdv_g<\infty$$ and $$\int_M|d\phi|^mdv_g<\epsilon$$ for $2<p<m$ and $\epsilon>0$ (depending on $p,A$ and the geometry of \(M\)) sufficiently small, then $\phi$ must be harmonic.
\end{thm}

The central idea in the proof of Theorem \ref{Bra}
is the use of an \emph{Euclidean type Sobolev inequality} that allows to control the curvature term
in the biharmonic map equation. However, in order for this inequality to hold one has to make stronger
assumptions on the domain manifold \(M\) as in Theorem \ref{Bra}, which we will correct below.

We say that a complete non-compact Riemannian manifold of infinite volume admits
an \emph{Euclidean type Sobolev inequality} if the following inequality holds (assuming \(m=\dim M\geq 3\))
\begin{align}
\label{sobolev-inequality}
(\int_M|u|^{2m/(m-2)}dv_g)^\frac{m-2}{m}\leq C_{sob}^M\int_M|\nabla u|^2dv_g
\end{align}
for all \(u\in W^{1,2}(M)\) with compact support,
where \(C_{sob}^M\) is a positive constant that depends on the geometry of \(M\).
Such an inequality holds in \(\R^m\) and is well-known as \emph{Gagliardo-Nirenberg inequality} in this case.

One way of ensuring that \eqref{sobolev-inequality} holds is the following:
If \((M,g)\) is a complete, non-compact Riemannian manifold of dimension \(m\)
with nonnegative Ricci curvature, and if for some point \(x\in M\)
\begin{align*}
\lim_{R\to\infty}\frac{vol_g(B_R(x))}{\omega_mR^m}>0
\end{align*}
holds, then \eqref{sobolev-inequality} holds true, see \cite{Sh}.
Here, \(\omega_m\) denotes the volume of the unit ball in \(\R^m\).
For further geometric conditions ensuring that \eqref{sobolev-inequality} holds
we refer to \cite[Section 3.7]{He}.

In this article we will correct the assumptions that are needed for Theorem \ref{Bra} to hold
and extend it to the case of $p=2$, which is a more natural integrability condition.
Motivated by these aspects, we actually can prove the following result:
\begin{thm}\label{main1}
Suppose that $(M,g)$ is a complete, connected non-compact Riemannian manifold of dimension \(m=\dim M\geq 3\) with infinite volume that admits
an Euclidean type Sobolev inequality of the form \eqref{sobolev-inequality}.
Moreover, suppose that \((N,h)\) is another Riemannian manifold
whose sectional curvatures satisfy $K^N\leq A,$ where $A$ is a positive constant.
Let $\phi: (M,g)\to (N,h)$ be a smooth biharmonic map.
If $$\int_M|\tau(\phi)|^pdv_g<\infty$$ and $$\int_M|d\phi|^mdv_g<\epsilon$$
for $p>1$ and $\epsilon>0$ (depending on $p,A$ and the geometry of \(M\)) sufficiently small, then $\phi$ must be harmonic.
\end{thm}

Similar ideas have been used to derive Liouville type results for \(p\)-harmonic maps in \cite{NT}, see also  \cite{zc} for a more general result.
In the proof we choose a test function of the form $(|\tau(\phi)|^2+\delta)^\frac{p-2}{2}\tau(\phi), (p>1, \delta>0)$ to avoid problems that may be caused by the zero points of $\tau(\phi)$.
When we take the limit $\delta\to 0$, we also need to be careful about the set of zero points of $\tau(\phi)$, and a delicate analysis is given.
For details please see the proof in section 2.

Moreover, we can get the following Liouville type result.
\begin{thm}\label{main2}
Suppose that $(M,g)$ is a complete, connected non-compact Riemannian manifold of \(m=\dim M\geq 3\) with nonnegative Ricci curvature that admits
an Euclidean type Sobolev inequality of the form \eqref{sobolev-inequality}.
Moreover, suppose that \((N,h)\) is another Riemannian manifold
whose sectional curvatures satisfy $K^N\leq A,$ where $A$ is a positive constant.
Let $\phi: (M,g)\to (N,h)$ be a smooth biharmonic map.
If $$\int_M|\tau(\phi)|^pdv_g<\infty$$ and $$\int_M|d\phi|^mdv_g<\epsilon$$
for $p>1$ and $\epsilon>0$ (depending on $p,A$ and the geometry of \(M\)) sufficiently small, then $\phi$ is a constant map.
\end{thm}
Note that due to a classical result of Calabi and Yau \cite[Theorem 7]{Yau} a complete non-compact Riemannian manifold
with nonnegative Ricci curvature has infinite volume.

\begin{rem}
Due to Theorem \ref{main1} we only need to prove that harmonic maps satisfying the assumption of Theorem \ref{main2} are constant maps.
Such a result was proven in \cite{NT} and thus Theorem \ref{main2} is a corollary of Theorem 1.5 in \cite{NT}. Conversely, Theorem \ref{main2} generalizes related Liouville type results for harmonic maps in \cite{NT}.
\end{rem}
\quad \\
\textbf{Organization:} Theorem \ref{main1} is proved in section 2.
In section 3 we apply Theorem \ref{main1} to get several nonexistence results for proper biharmonic submersions.

\section{Proof of the main result}
In this section we will prove Theorem \ref{main1}.

Assume that $x_0\in M$. We choose a cutoff function $0\leq\eta\leq1$ on $M$ that satisfies
\begin{equation}\label{flow1}
\left\{\begin{array}{rcl}
\eta(x)&=&1, \quad \forall \ x\in B_R(x_0), \\
\eta(x)&=&0,\quad \forall \ x\in M\setminus B_{2R}(x_0),\\
|\nabla\eta(x)|&\leq& \frac{C}{R}, \quad \forall \  x \in M.
\end{array}\right.
\end{equation}

\begin{lem}\label{lem1}
Let $\phi:(M,g)\to (N,h)$ be a smooth biharmonic map and assume that the sectional curvatures of $N$ satisfy $K^N\leq A.$
Let $\delta$ be a positive constant. Then the following inequalities hold.
\\(1) If $1<p<2$, we have
\begin{eqnarray}\label{ine1}
&&(1-\frac{p-1}{2})\int_M\eta^2(|\tau(\phi)|^2+\delta)^\frac{p-4}{2}|\bar{\nabla}\tau(\phi)|^2|\tau(\phi)|^2dv_g \nonumber
\\&\leq& A\int_M\eta^2(|\tau(\phi)|^2+\delta)^\frac{p}{2}|d\phi|^2dv_g+\frac{C}{R^2}\int_{B_{2R}(x_0)}(|\tau(\phi)|^2+\delta)^\frac{p}{2}dv_g\nonumber
\\&-&(p-2)\int_M\eta^2(|\tau(\phi)|^2+\delta)^\frac{p-4}{2}|\bar{\nabla}\tau(\phi)|^2|\tau(\phi)|^2dv_g;
\end{eqnarray}
\\(2) If $p\geq2$, we have
\begin{eqnarray}\label{ine1'}
&&\frac{1}{2}\int_M\eta^2(|\tau(\phi)|^2+\delta)^\frac{p-4}{2}|\bar{\nabla}\tau(\phi)|^2|\tau(\phi)|^2dv_g \nonumber
\\&\leq& A\int_M\eta^2(|\tau(\phi)|^2+\delta)^\frac{p}{2}|d\phi|^2dv_g+\frac{C}{R^2}\int_{B_{2R}(x_0)}(|\tau(\phi)|^2+\delta)^\frac{p}{2}dv_g.
\end{eqnarray}
\end{lem}
\proof  Multiplying the biharmonic map equation by $\eta^2(|\tau(\phi)|^2+\delta)^\frac{p-2}{2}\tau(\phi)$ we get
$$\eta^2(|\tau(\phi)|^2+\delta)^\frac{p-2}{2}\langle\Delta^\phi\tau(\phi), \tau(\phi)\rangle=-\eta^2(|\tau(\phi)|^2+\delta)^\frac{p-2}{2}\sum_{i=1}^mR^N(\tau(\phi),d\phi(e_i),\tau(\phi),d\phi(e_i)).$$

Integrating over $M$ and using integration by parts we get
\begin{eqnarray}\label{inem1}
&&\int_M\eta^2(|\tau(\phi)|^2+\delta)^\frac{p-2}{2}\langle\Delta^\phi\tau(\phi), \tau(\phi)\rangle dv_g\nonumber
\\&=&-2\int_M(|\tau(\phi)|^2+\delta)^\frac{p-2}{2}\langle\bar{\nabla}\tau(\phi), \tau(\phi)\rangle\eta\nabla\eta dv_g\nonumber
\\&-&(p-2)\int_M\eta^2|\langle\bar{\nabla}\tau(\phi), \tau(\phi)\rangle|^2(|\tau(\phi)|^2+\delta)^\frac{p-4}{2}dv_g\nonumber
\\&-&\int_M\eta^2(|\tau(\phi)|^2+\delta)^\frac{p-2}{2}|\bar{\nabla}\tau(\phi)|^2dv_g\nonumber
\\&\leq&-2\int_M(|\tau(\phi)|^2+\delta)^\frac{p-2}{2}\langle\bar{\nabla}\tau(\phi), \tau(\phi)\rangle\eta\nabla\eta dv_g
\\&-&(p-2)\int_M\eta^2|\langle\bar{\nabla}\tau(\phi), \tau(\phi)\rangle|^2(|\tau(\phi)|^2+\delta)^\frac{p-4}{2}dv_g\nonumber
\\&-&\int_M\eta^2(|\tau(\phi)|^2+\delta)^\frac{p-4}{2}|\bar{\nabla}\tau(\phi)|^2|\tau(\phi)|^2dv_g.\nonumber
\end{eqnarray}

Therefore when $1<p<2$ we have

\begin{eqnarray*}
&&\int_M\eta^2(|\tau(\phi)|^2+\delta)^\frac{p-2}{2}\langle\Delta^\phi\tau(\phi), \tau(\phi)\rangle dv_g
\\&\leq&(\frac{p-1}{2}-1)\int_M\eta^2(|\tau(\phi)|^2+\delta)^\frac{p-4}{2}|\bar{\nabla}\tau(\phi)|^2|\tau(\phi)|^2dv_g
\\&+&\frac{2}{p-1}\int_M(|\tau(\phi)|^2+\delta)
^\frac{p}{2}|\nabla\eta|^2dv_g
\\&-&(p-2)\int_M\eta^2|\langle\bar{\nabla}\tau(\phi), \tau(\phi)\rangle|^2(|\tau(\phi)|^2+\delta)^\frac{p-4}{2}dv_g
\\&\leq& (\frac{p-1}{2}-1)\int_M\eta^2(|\tau(\phi)|^2+\delta)^\frac{p-4}{2}|\bar{\nabla}\tau(\phi)|^2|\tau(\phi)|^2dv_g
\\&+&\frac{C}{R^2}\int_{B_{2R(x_0)}}(|\tau(\phi)|^2+\delta)
^\frac{p}{2}dv_g
\\&-&(p-2)\int_M\eta^2|\langle\bar{\nabla}\tau(\phi), \tau(\phi)\rangle|^2(|\tau(\phi)|^2+\delta)^\frac{p-4}{2}dv_g
\\&\leq& (\frac{p-1}{2}-1)\int_M\eta^2(|\tau(\phi)|^2+\delta)^\frac{p-4}{2}|\bar{\nabla}\tau(\phi)|^2|\tau(\phi)|^2dv_g
\\&+&\frac{C}{R^2}\int_{B_{2R(x_0)}}(|\tau(\phi)|^2+\delta)
^\frac{p}{2}dv_g
\\&-&(p-2)\int_M\eta^2|\bar{\nabla}\tau(\phi)|^2|\tau(\phi)|^2(|\tau(\phi)|^2+\delta)^\frac{p-4}{2}dv_g,
\end{eqnarray*}
where in the last inequality we used $1<p<2$ and
$$|\langle\bar{\nabla}\tau(\phi), \tau(\phi)\rangle|^2\leq|\bar{\nabla}\tau(\phi)|^2|\tau(\phi)|^2.$$
Therefore we find
\begin{eqnarray*}
&&(1-\frac{p-1}{2})\int_M\eta^2(|\tau(\phi)|^2+\delta)^\frac{p-4}{2}|\bar{\nabla}\tau(\phi)|^2|\tau(\phi)|^2dv_g \nonumber
\\&\leq& \int_M\eta^2(|\tau(\phi)|^2+\delta)^\frac{p-2}{2}\sum_{i=1}^mR^N(\tau(\phi),d\phi(e_i),\tau(\phi),d\phi(e_i))dv_g
\\&+&\frac{C}{R^2}\int_{B_{2R(x_0)}}(|\tau(\phi)|^2+\delta)
^\frac{p}{2}dv_g\nonumber
\\&-&(p-2)\int_M\eta^2(|\tau(\phi)|^2+\delta)^\frac{p-4}{2}|\bar{\nabla}\tau(\phi)|^2|\tau(\phi)|^2dv_g
\\&\leq& A\int_M\eta^2(|\tau(\phi)|^2+\delta)^\frac{p-2}{2}|\tau(\phi)|^2|d\phi|^2dv_g+\frac{C}{R^2}\int_{B_{2R(x_0)}}(|\tau(\phi)|^2+\delta)
^\frac{p}{2}dv_g\nonumber
\\&-&(p-2)\int_M\eta^2(|\tau(\phi)|^2+\delta)^\frac{p-4}{2}|\bar{\nabla}\tau(\phi)|^2|\tau(\phi)|^2dv_g
\\&\leq& A\int_M\eta^2(|\tau(\phi)|^2+\delta)^\frac{p}{2}|d\phi|^2dv_g+\frac{C}{R^2}\int_{B_{2R(x_0)}}(|\tau(\phi)|^2+\delta)
^\frac{p}{2}dv_g\nonumber
\\&-&(p-2)\int_M\eta^2(|\tau(\phi)|^2+\delta)^\frac{p-4}{2}|\bar{\nabla}\tau(\phi)|^2|\tau(\phi)|^2dv_g,
\end{eqnarray*}
which proves the first claim.\\

When $p\geq2$ equation \eqref{inem1} gives
\begin{eqnarray*}
&&\int_M\eta^2(|\tau(\phi)|^2+\delta)^\frac{p-2}{2}\langle\Delta^\phi\tau(\phi), \tau(\phi)\rangle dv_g
\\&\leq&-2\int_M(|\tau(\phi)|^2+\delta)^\frac{p-2}{2}\langle\bar{\nabla}\tau(\phi), \tau(\phi)\rangle\eta\nabla\eta dv_g
\\&-&\int_M\eta^2(|\tau(\phi)|^2+\delta)^\frac{p-4}{2}|\bar{\nabla}\tau(\phi)|^2|\tau(\phi)|^2dv_g
\\&\leq&-\frac{1}{2}\int_M\eta^2(|\tau(\phi)|^2+\delta)^\frac{p-4}{2}|\bar{\nabla}\tau(\phi)|^2|\tau(\phi)|^2dv_g+2\int_M(|\tau(\phi)|^2+\delta)
^\frac{p}{2}|\nabla\eta|^2dv_g.
\end{eqnarray*}

Therefore we have
\begin{eqnarray*}
&&\frac{1}{2}\int_M\eta^2(|\tau(\phi)|^2+\delta)^\frac{p-4}{2}|\bar{\nabla}\tau(\phi)|^2|\tau(\phi)|^2dv_g
\\&\leq&\int_M\eta^2(|\tau(\phi)|^2+\delta)^\frac{p-2}{2}\sum_{i=1}^mR^N(\tau(\phi),d\phi(e_i),\tau(\phi),d\phi(e_i))dv_g+\frac{C}{R^2}\int_{B_{2R(x_0)}}(|\tau(\phi)|^2+\delta)
^\frac{p}{2}dv_g
\\&\leq& A\int_M\eta^2(|\tau(\phi)|^2+\delta)^\frac{p-2}{2}|\tau(\phi)|^2|d\phi|^2dv_g+\frac{C}{R^2}\int_{B_{2R(x_0)}}(|\tau(\phi)|^2+\delta)
^\frac{p}{2}dv_g
\\&\leq& A\int_M\eta^2(|\tau(\phi)|^2+\delta)^\frac{p}{2}|d\phi|^2dv_g+\frac{C}{R^2}\int_{B_{2R(x_0)}}(|\tau(\phi)|^2+\delta)
^\frac{p}{2}dv_g.
\end{eqnarray*}

This completes the proof of Lemma \ref{lem1}.
\endproof

In the following we will estimate the term $$A\int_M\eta^2(|\tau(\phi)|^2+\delta)^\frac{p}{2}|d\phi|^2dv_g.$$
\begin{lem}\label{lem2}
Assume that $(M,g)$ satisfies the assumptions of Theorem \ref{main1}. Then the following inequality holds
\begin{eqnarray}\label{ine2}
&&\int_M\eta^2(|\tau(\phi)|^2+\delta)^\frac{p}{2}|d\phi|^2dv_g\nonumber
\\&\leq & C(\int_M|d\phi|^mdv_g)^{\frac{2}{m}}\times
\\&&(\frac{1}{R^2}\int_{B_{2R(x_0)}}(|\tau(\phi)|^2+\delta)^\frac{p}{2}dv_g\nonumber
+\int_M\eta^2(|\tau(\phi)|^2+\delta)^\frac{p-4}{2}|\bar{\nabla}\tau(\phi)|^2|\tau(\phi)|^2dv_g),
\end{eqnarray}
where $C$ is a constant depending on $p,A$ and the geometry of $M$.
\end{lem}

\proof Set $f=(|\tau(\phi)|^2+\delta)^\frac{p}{4}$, then we have
$$\int_M\eta^2(|\tau(\phi)|^2+\delta)^\frac{p}{2}|d\phi|^2dv_g=\int_M\eta^2f^2|d\phi|^2dv_g.$$ Then by H\"older's inequality we get
$$\int_M\eta^2f^2|d\phi|^2dv_g\leq (\int_M(\eta f)^\frac{2m}{m-2}dv_g)^\frac{m-2}{m}(\int_M|d\phi|^mdv_g)^\frac{2}{m}.$$

Applying \eqref{sobolev-inequality} to $u=\eta f$ we get
$$(\int_M(\eta f)^\frac{2m}{m-2}dv_g)^\frac{m-2}{m}\leq C^M_{sob}\int_M|d(\eta f)|^2dv_g,$$
which leads to
\begin{eqnarray}\label{ine7}
\int_M\eta^2f^2|d\phi|^2dv_g\leq 2C^M_{sob}(\int_M|d\phi|^mdv_g)^{\frac{2}{m}}(\int_M|d\eta|^2f^2dv_g+\int_M\eta^2|df|^2dv_g).
\end{eqnarray}
Note that $f=(|\tau(\phi)|^2+\delta)^\frac{p}{4}$ and $$|df|^2=\frac{p^2}{4}(|\tau(\phi)|^2+\delta)^\frac{p-4}{2}|\langle\bar{\nabla}\tau(\phi), \tau(\phi)\rangle|^2\leq \frac{p^2}{4}(|\tau(\phi)|^2+\delta)^\frac{p-4}{2}|\bar{\nabla}\tau(\phi)|^2|\tau(\phi)|^2.$$
This completes the proof of Lemma \ref{lem2}.
\endproof

When $1<p<2$, due to Lemmas \ref{lem1}, \ref{lem2}, we see that by choosing $\epsilon$ sufficiently small such that $AC\epsilon^\frac{2}{m}\leq\frac{p-1}{4},$ we have
\begin{eqnarray}\label{ine3}
\frac{p-1}{4}\int_M\eta^2(|\tau(\phi)|^2+\delta)^\frac{p-4}{2}|\bar{\nabla}\tau(\phi)|^2|\tau(\phi)|^2dv_g
\leq\frac{C}{R^2}\int_{B_{2R(x_0)}}(|\tau(\phi)|^2+\delta)^\frac{p}{2}dv_g,
\end{eqnarray}
where $C$ is a constant depending on $p,A$ and the geometry of $M$.

Now, set $M_1:=\{x\in M| \tau(\phi)(x)=0\}$, and $M_2=M\setminus M_1.$

If $M_2$ is an empty set, then we are done. Hence we assume that $M_2$ is nonempty and we will get a contradiction below. 

Note that since $\phi$ is smooth, $M_2$ is an open set.

From \eqref{ine3} we have
\begin{eqnarray}
\frac{p-1}{4}\int_{M_2}\eta^2(|\tau(\phi)|^2+\delta)^\frac{p-4}{2}|\bar{\nabla}\tau(\phi)|^2|\tau(\phi)|^2dv_g
\leq\frac{C}{R^2}\int_{B_{2R(x_0)}}(|\tau(\phi)|^2+\delta)^\frac{p}{2}dv_g.
\end{eqnarray}
Letting $\delta\to 0$ we get
$$\frac{p-1}{4}\int_{M_2}\eta^2|\tau(\phi)|^{p-2}|\bar{\nabla}\tau(\phi)|^2dv_g\leq\frac{C}{R^2}\int_{B_{2R(x_0)}}|\tau(\phi)|^pdv_g\leq
\frac{C}{R^2}\int_{M}|\tau(\phi)|^pdv_g.$$
Letting $R\to \infty$ we get
$$\frac{p-1}{4}\int_{M_2}|\tau(\phi)|^{p-2}|\bar{\nabla}\tau(\phi)|^2dv_g=0.$$
When $p\geq2$ by a similar discussion	 we can prove that
$$\frac{1}{4}\int_{M_2}|\tau(\phi)|^{p-2}|\bar{\nabla}\tau(\phi)|^2dv_g=0.$$
Therefore we have that $\bar{\nabla}\tau(\phi)=0$ everywhere in $M_2$ and hence $M_2$ is an \textbf{open} and \textbf{closed} nonempty set, thus $M_2=M$ (as we assume that $M$ is a connected manifold)
and $|\tau(\phi)|\equiv c$ for some constant $c\neq 0$. Thus $Vol(M)<\infty$ by $\int_Mc^pdv_g<\infty$.

In the following we will need Gaffney's theorem \cite{Ga}, stated below:
\begin{thm}[Gaffney]
Let $(M, g)$ be a complete Riemannian manifold. If a $C^1$ 1-form $\omega$ satisfies
that $\int_M|\omega|dv_g<\infty \ and \ \int_M|\delta\omega| dv_g<\infty,$ or equivalently, a $C^1$ vector field $X$ defined by
$\omega(Y) = \langle X, Y \rangle, (\forall Y \in TM)$ satisfies that $\int_M|X|dv_g<\infty \ and \ \int_M|div X|dv_g<\infty,$ then $$\int_M\delta\omega dv_g=\int_Mdiv Xdv_g=0.$$
\end{thm}

Define a l-form on $M$ by
$$\omega(X):=\langle d\phi(X),\tau(\phi)\rangle,~(X\in TM).$$
Then
\begin{eqnarray*}
\int_M|\omega|dv_g&=&\int_M(\sum_{i=1}^m|\omega(e_i)|^2)^\frac{1}{2}dv_g
\\&\leq&\int_M|\tau(\phi)||d\phi|dv_g
\\&\leq&c Vol(M)^{1-\frac{1}{m}}(\int_M|d\phi|^mdv_g)^\frac{1}{m}
\\&<&\infty.
\end{eqnarray*}
In addition, we calculate $-\delta\omega=\sum_{i=1}^m(\nabla_{e_i}\omega)(e_i)$:
\begin{eqnarray*}
-\delta\omega&=&\sum_{i=1}^m\nabla_{e_i}(\omega(e_i))-\omega(\nabla_{e_i}e_i)
\\&=&\sum_{i=1}^m\{\langle\bar{\nabla}_{e_i}d\phi(e_i),\tau(\phi)\rangle
-\langle d\phi(\nabla_{e_i}e_i),\tau(\phi)\rangle\}
\\&=&\sum_{i=1}^m\langle \bar{\nabla}_{e_i}d\phi(e_i)-d\phi(\nabla_{e_i}e_i),\tau(\phi)\rangle
\\&=&|\tau(\phi)|^2,
\end{eqnarray*}
where in the second equality we used $\bar{\nabla}\tau(\phi)=0$. Therefore $$\int_M|\delta\omega|dv_g=c^2Vol(M)<\infty.$$
Now by Gaffney's theorem and the above equality we have that
$$0=\int_M(-\delta\omega)dv_g=\int_M|\tau(\phi)|^2dv_g=c^2Vol(M),$$
which implies that $c=0$, a contradiction. Therefore we must have $M_1=M$, i.e. $\phi$ is a harmonic map. This completes the proof of Theorem \ref{main1}.
\endproof
\section{Applications to biharmonic submersions}
In this section we give some applications of our result to biharmonic submersions.

First we recall some definitions \cite{BW}.

Assume that $\phi: (M, g)\to (N, h)$ is a smooth map between Riemannian manifolds and $x\in M$. Then $\phi$ is called {\bf horizontally weakly conformal} if either

(i) $d\phi_x=0$, or

(ii) $d\phi_x$ maps the horizontal space $\rm \mathcal{H}_x=\{Ker~d\phi_x\}^\bot$ conformally \textbf{onto} $T_{\phi(x)}N$, i.e.
$$h(d\phi_x(X), d\phi_x(Y))=\lambda^2 g(X, Y), (X, Y\in \mathcal{H}_x),$$
for some $\lambda=\lambda(x)>0,$ called the {\bf dilation} of $\phi$ at $x$.

A map $\phi$ is called {\bf horizontally weakly conformal} or {\bf semiconformal} on $M$ if it is horizontally weakly conformal at every point of $M$. Furthermore, if $\phi$ has no critical points, then we call it a {\bf horizontally conformal submersion}: In this case the dilation $\lambda:M \to (0,\infty)$ is a smooth function. Note that if $\phi: (M, g)\to (N, h)$ is a horizontally weakly conformal map and $\dim M<\dim N$, then $\phi$ is a constant map.

If for every harmonic function $f: V\to \mathbb{R}$ defined on an open subset $V$ of $N$ with $\phi^{-1}(V)$ nonempty, the composition $f\circ\phi$ is harmonic on $\phi^{-1}(V)$, then $\phi$ is called a {\bf harmonic morphism}. Harmonic morphisms are characterized as follows (cf. \cite{Fu, Is}).
\begin{thm}[\cite{Fu, Is}]\label{thm4}
A smooth map $\phi: (M, g)\to (N, h)$ between Riemannian manifolds is a harmonic morphism if and only if $\phi$ is both harmonic and horizontally weakly conformal.
\end{thm}

When $\phi:(M^m, g)\to (N^n, h),(m>n\geq2)$ is a horizontally conformal submersion, the tension field is given by
\begin{eqnarray}\label{eq5}
\tau(\phi)=\frac{n-2}{2}\lambda^2d\phi(grad_\mathcal{H}(\frac{1}{\lambda^2}))
-(m-n)d\phi(\hat{H}),
\end{eqnarray}
where $grad_\mathcal{H}(\frac{1}{\lambda^2})$ is the horizontal component of $\rm grad(\frac{1}{\lambda^2})$, and $\hat{H}$ is the {\bf mean curvature} of the fibres given by the trace
$$\hat{H}=\frac{1}{m-n}\sum_{i=n+1}^m\mathcal{H}(\nabla_{e_i}e_i).$$
Here, $\{e_i, i=1,...,m\}$ is a local orthonormal frame field on $M$ such that $\{e_{i}, i=1,...,n\}$ belongs to $\mathcal{H}_x$ and $\{e_{j}, j=n+1,...,m \}$ belongs to $\mathcal{V}_x$ at each point $x\in M$, where $T_xM=\mathcal{H}_x\oplus \mathcal{V}_x$.

Nakauchi et al. \cite{NUG}, Maeta \cite{Ma} and Luo \cite{Luo2} applied their nonexistence result for biharmonic maps to get conditions for which biharmonic submersions are harmonic morphisms.
Here, we give another such result by using Theorem \ref{main1}. We have
\begin{pro}
Let $\phi:(M^m, g)\to (N^n, h), (m>n\geq2)$ be a biharmonic horizontally conformal submersion from a complete, connected non-compact Riemannian manifold $(M,g)$ with infinite volume,
that admits an Euclidean Sobolev type inequality of the form \eqref{sobolev-inequality},
into a Riemannian manifold $(N, h)$ with sectional curvatures $K^N\leq A$ and $p$ a real constant satisfying $1<p<\infty$. If
 $$ \int_M\lambda^p|\frac{n-2}{2}\lambda^2grad_\mathcal{H}(\frac{1}{\lambda^2})
-(m-n)\hat{H}|_g^pdv_g<\infty,$$
and $$\int_M\lambda^mdv_g<\epsilon$$ for sufficiently small $\epsilon>0$ (depending on $p,A$ and the geometry of \(M\)), then $\phi$ is a harmonic morphism.
\end{pro}
\proof By (\ref{eq5}) we have,
$$\int_M|\tau(\phi)|_h^pdv_g=\int_M\lambda^p|\frac{n-2}{2}\lambda^2grad_\mathcal{H}(\frac{1}
{\lambda^2})
-(m-n)\hat{H}|_g^pdv_g<\infty,$$
and since $|d\phi(x)|^2=n\lambda^2(x)$, we get that $\phi$ is harmonic by Theorem \ref{main1}. Since $\phi$ is also a horizontally conformal submersion, $\phi$ is a harmonic morphism by Theorem \ref{thm4}. \endproof

In particular, if $\dim N=2$, we have
\begin{cor}
Let $\phi:(M^m, g)\to (N^2, h)$ be a biharmonic horizontally conformal submersion from a complete, connected non-compact Riemannian manifold $(M, g)$ with infinite volume,
that admits an Euclidean Sobolev type inequality of the form \eqref{sobolev-inequality},
into a Riemannian surface $(N, h)$ with Gauss curvature bounded from above and $p$ a real constant satisfying $1<p<\infty$. If
$$\int_M\lambda^p|\hat{H}|_g^pdv_g<\infty,$$
and $$\int_M\lambda^mdv_g<\epsilon$$ for sufficiently small $\epsilon>0$ (depending on $p,A$ and the geometry of \(M\)), then $\phi$ is a harmonic morphism.
\end{cor}

\quad\\

\textbf{Acknowledgements:}
The first named author gratefully acknowledges the support of the Austrian Science Fund (FWF)
through the project P30749-N35 ``Geometric variational problems from string theory''.
The second named author is supported by the NSF of China(No.11501421, No.11771339). Part of the work was finished when the second named author
was a visiting scholar at Tsinghua University. He would like to express his gratitude to
Professor Yuxiang Li and Professor Hui Ma for their invitation and to Tsinghua University for their hospitality. The second named author also would like to thank Professor Ye Lin Ou for his interest in this work and discussion. 

{}
\vspace{1cm}\sc
Volker Branding

Faculty of Mathematics,

University of Vienna, Oskar-Morgenstern-Platz 1, 1090 Vienna, Austria

{\tt volker.branding@univie.ac.at}

\vspace{1cm}\sc

Yong Luo

School of mathematics and statistics,

Wuhan university, Wuhan 430072, China

{\tt yongluo@whu.edu.cn}

\vspace{1cm}\sc

\end{document}